\documentclass[12pt]{amsart}
\usepackage{amsthm, amsmath, amssymb, amsfonts, url}
\usepackage[margin = 1in]{geometry}

\usepackage{hyperref}

\hypersetup{
    pdfstartview={FitH}    
}

\newtheorem*{claim*}{Claim}
\newtheorem*{fact*}{Fact}
\newtheorem{theorem}{Theorem}
\newtheorem*{theorem*}{Theorem}

\newtheorem*{lemma*}{Lemma}
\newtheorem{corollary}[theorem]{Corollary}
\newtheorem{conjecture}[theorem]{Conjecture}

\theoremstyle{definition}

\theoremstyle{remark}
\newtheorem*{remark}{Remark}

\def\c{\mathcal}

\DeclareMathOperator{\Hom}{Hom}

\newcommand{\abs}[1]{\left\lvert#1\right\rvert}

\def\({\left(}
\def\){\right)}

\def\x{\times}

\def\<={\Leftarrow}
\def\=>{\Rightarrow}

\title{The Number of Independent Sets in a Regular Graph}
\author{Yufei Zhao}
\address{Department of Mathematics, Massachusetts Institute of Technology, Cambridge, MA 02139}
\email{yufeiz@mit.edu}
\date{August 2, 2009}

\begin{document}

\begin{abstract}
We show that the number of independent sets in an $N$-vertex, $d$-regular graph is at most $(2^{d+1} - 1 )^{N/2d}$, where the bound is sharp for a disjoint union of complete $d$-regular bipartite graphs. This settles a conjecture of Alon in 1991 and Kahn in 2001. Kahn proved the bound when the graph is assumed to be bipartite. We give a short proof that reduces the general case to the bipartite case. Our method also works for a weighted generalization, i.e., an upper bound for the independence polynomial of a regular graph.
\end{abstract}

\maketitle

\section{Introduction}

Let $G$ be a (simple, finite, undirected) graph. An \emph{independent set} is a subset of the vertices with no two adjacent. Let $\c I(G)$ denote the collection of independent sets of $G$ and let $i(G)$ be its cardinality. We would like to determine an upper bound for $i(G)$ when $G$ is an $N$-vertex, $d$-regular graph. This problem has received some attention due to its applications in combinatorial group theory \cite{Alon91} and statistical mechanics \cite{Kahn}. We begin by giving a brief overview of the history of the problem (also see \cite{Galvin}).

In the 1988 Number Theory Conference at Banff, in an effort to resolve the Cameron-Erd\H{o}s conjecture \cite{CE, CE-Green, CE-Sap} on the number of sum-free sets, A.~Granville conjectured that $i(G) \leq 2^{(1/2 + o(1))N}$ for $N$-vertex, $d$-regular graphs $G$, where $o(1) \to 0$ as $d \to \infty$. In 1991, Alon \cite{Alon91} resolved Granville's conjecture by proving that 
\begin{equation} \label{eq:Alon}
	i(G) \leq 2^{(1/2 + O(d^{-0.1}))N}
\end{equation}
and applied the result to counting sum-free sets in groups. Alon speculated that perhaps $i(G)$ is maximized (for $N$ divisible by $2d$) when $G$ is a disjoint union of $N/2d$ complete $d$-regular bipartite graphs $K_{d,d}$. This conjecture was later formally stated by Kahn \cite{Kahn} in 2001:

\begin{conjecture}[Alon and Kahn] \label{conj:Kahn}
For any $N$-vertex, $d$-regular graph $G$,
\begin{equation} \label{eq:Kahn-conj}
	i(G) \leq i\(K_{d,d}\)^{N/2d} = \(2^{d+1} - 1 \)^{N/2d}
\end{equation}
\end{conjecture}

Here we prove this conjecture as a special case of our main result, Theorem \ref{thm:main} below. Kahn used entropy methods to prove Conjecture \ref{conj:Kahn} when $G$ is bipartite and applied it to prove some results about the hard-core lattice gas model in statistical mechanics. Several improvements to \eqref{eq:Alon} for general $G$ were later given by Sapozhenko \cite{Sapo}, Kahn (see \cite{MT}), and Galvin \cite{Galvin}. The best previous bound, given recently by Galvin \cite{Galvin}, is
\[
	i(G) \leq \exp_2 \(\frac{N}{2} \( 1 + \frac{1}{d} +  \frac{O(1)}{d} \sqrt{\frac{\log d}{d}}\) \).
\]
For comparison, the bound in \eqref{eq:Kahn-conj}, which is tight whenever $N$ is divisible by $2d$, has the asymptotic expansion
\[
	\(2^{d+1} - 1 \)^{N/2d} = \exp_2 \(\frac{N}{2} \( 1 + \frac 1 d - \frac{1+o(1)}{(2\ln 2) d 2^d} \) \),
\]
where $o(1) \to 0$ as $d \to \infty$. The bounds of \cite{Galvin} and \cite{Sapo} only work well for large $d$. On the other hand, recent work by Galvin and the author \cite{GZ} showed through computational verification that Conjecture \ref{conj:Kahn} is true for $d \leq 5$.

We also consider a weighted generalization of Conjecture \ref{conj:Kahn}. The \emph{independence polynomial} of $G$, introduced by Gutman and Harary \cite{GutHar}, is defined as
\[
	P(\lambda, G) = \sum_{I \in \c I(G)} \lambda^{\abs{I}}.
\]
This is also the partition function for the \emph{hard-core model} on $G$ with \emph{activity} (or \emph{fugacity}) $\lambda$ (see \cite{BergSteif, Haggstrom} for background). Galvin and Tetali \cite{GT} applied Kahn's techniques to weighted graph homomorphisms. In particular, a generalization of Conjecture \ref{conj:Kahn} was obtained for $N$-vertex, $d$-regular bipartite $G$:
\begin{equation} \label{eq:hardcore}
	P(\lambda, G) \leq \(2 (1+\lambda)^d - 1 \)^{N/2d}
\end{equation}
for all $\lambda \geq 0$, where equality is once again attained when $G$ is a disjoint union of $K_{d,d}$'s. It has been conjectured in \cite{Galvin} and \cite{GT} that \eqref{eq:hardcore} holds also for non-bipartite graphs. In \cite{CGT} and \cite{Galvin}, upper bounds were given for $P(\lambda, G)$ for general $G$. In this paper, we settle the conjecture by proving \eqref{eq:hardcore} for all graphs $G$.

\begin{theorem} \label{thm:main}
For any $N$-vertex, $d$-regular graph $G$, and any $\lambda \geq 0$,
\[
	P(\lambda, G) \leq P\(\lambda, K_{d,d} \)^{N/2d} = \(2 (1+\lambda)^d - 1 \)^{N/2d}.
\]
\end{theorem}

When $\lambda = 1$, we obtain Conjecture \ref{conj:Kahn} as a corollary.

\begin{corollary} \label{cor:main}
For any $N$-vertex, $d$-regular graph $G$,
\[
	i(G) \leq i\(K_{d,d}\)^{N/2d} = \(2^{d+1} - 1 \)^{N/2d}.
\]
\end{corollary}

In Section \ref{sec:proof} we derive Theorem \ref{thm:main} from its specialization to bipartite graphs (which was known earlier). Section \ref{sec:discussion} discusses some further questions.

\section{Proof of Theorem \ref{thm:main}} \label{sec:proof}

In this section, we prove Theorem \ref{thm:main} by reducing to the bipartite case, which was proven in \cite{GT} (see \cite{Kahn} for the non-weighted case).

Let $V(G)$ denote the set of vertices of $G$. For $A, B \subset V(G)$, say that $A$ is \emph{independent from} $B$ if $G$ contains no edge of the form $ab$, where $a \in A$ and $b \in B$. Also let $G[A]$ denote the subgraph of $G$ induced by $A$. 

Let $\c J(G)$ denote the set of pairs $(A,B)$ of subsets of vertices of $G$ such that $A$ is independent from $B$ and $G[A \cup B]$ is bipartite. For a pair $(A, B)$ of subsets of $V(G)$, define its \emph{size} to be $\abs{A} + \abs{B}$.

Here is the key lemma of our proof.

\begin{lemma*}
For any graph $G$ (not necessarily regular), there exists a size-preserving bijection between $\c I(G) \x \c I(G)$ and $\c J(G)$.
\end{lemma*}

\begin{proof}
For every $W \subset V(G)$ such that $G[W]$ is bipartite, fix a bipartition $W = W_1 \cup W_2$ so that $W_1$ and $W_2$ are both independent in $G$. Let $\c K(G)$ be the set of pairs $(A, B)$ of subsets of $V(G)$ such that $G[A\cup B]$ is bipartite. Note that $\c I(G) \x \c I(G) \subset K(G)$ and $\c J(G) \subset \c K(G)$. The former is true because if $A, B \in \c I(G)$ then $A \cup (B \setminus A)$ is already a bipartition.

We construct a size-preserving involution of $\c K(G)$ as follows. For any $(A, B) \in \c K$, let $W_1 \cup W_2$ be the chosen bipartition of $W = A \cup B$. Then let the involution send $(A, B)$ to
\[
	((A \cap W_1) \cup (B \cap W_2),  (A \cap W_2) \cup (B \cap W_1)).
\]
It is easy to check that this is a size-preserving involution of $\c K(G)$ that maps $\c I(G) \x \c I(G)$ to $\c J(G)$ and vice-versa. The lemma follows immediately.
\end{proof}


%

\begin{proof}[Proof of Theorem \ref{thm:main}]
Let $G \x K_2$ denote the bipartite double cover of $G$. This is the bipartite graph with vertices $(v, i)$, where $v \in V(G)$ and $i \in \{0, 1\}$, and an edge between $(v, 0)$ and $(w, 1)$ whenever $vw$ is an edge in $G$.
Then, independent sets in $G \x K_2$ correspond to pairs $(A,B)$ of subsets of $V(G)$ such that $A$ is independent from $B$. Hence for any $\lambda \geq 0$,
\begin{equation} \label{eq:1}
	P(\lambda, G \x K_2)
	= \hspace{-1em} \sum_{I \in \c I(G \x K_2)}  \hspace{-1em} \lambda^{\abs{I}}
	= \hspace{-1.5em} \sum_{\substack{A, B \subset V(G) \\ A \text{ indep.~from } B}} \hspace{-1.5em}  \lambda^{\abs{A} + \abs{B}}
	\geq \hspace{-1em} \sum_{(A,B) \in \c J(G)} \hspace{-1em} \lambda^{\abs{A} + \abs{B}}
	= \hspace{-1em} \sum_{A, B \in \c I(G)} \hspace{-1em} \lambda^{\abs{A} + \abs{B}}
	= P(\lambda, G)^2,
\end{equation}
where we used the lemma at the second to last equality. When $G$ is a $d$-regular graph, so is $G \x K_2$. Note that $G \x K_2$ is bipartite, and we know that the theorem is true for bipartite graphs, so
\begin{equation} \label{eq:2}
	P(\lambda, G \x K_2) \leq P(\lambda, K_{d,d} )^{N/d}.
\end{equation}
Combining \eqref{eq:1} and \eqref{eq:2} gives us
\[
	P(\lambda, G) \leq P(\lambda, K_{d,d} )^{N/2d},
\]
as desired.
\end{proof}

\begin{remark}
Kahn \cite{Kahn} also conjectured that, for any graph $G$ without isolated vertices
\begin{equation} \label{eq:edge-weighted}
	i(G) \leq \prod_{uv \in E(G)} \(2^{d(u)} + 2^{d(v)} - 1\)^{1/d(u)d(v)},
\end{equation}
where $E(G)$ denotes the set of edges of $G$ and $d(v)$ the degree of vertex $v$. Galvin \cite{GalvinPC} noted that the proof presented here can also be used to reduce \eqref{eq:edge-weighted} to the bipartite case. Recently Galvin and the author \cite{GZ} showed that \eqref{eq:edge-weighted} is true if the maximum degree of $G$ is at most 5.
\end{remark}

\section{Discussion and Further questions} \label{sec:discussion}

\subsection*{Non-entropy proof of the bipartite case?}
So far the only known proofs of Theorem \ref{thm:main} and Corollary \ref{cor:main} in the bipartite case use entropy methods. The entropy approach was first applied to this problem by Kahn \cite{Kahn}, and was subsequently applied to several variations of the problem. For instance, Kahn \cite{Kahn2} later showed that for an $N$-vertex, $d$-regular bipartite graph $G$ with vertex bipartition $\c O \cup \c E$, we have
\begin{equation} \label{eq:biweighted}
	\sum_{I \in \c I(G)} \mu^{\abs{I \cap \c O}} \lambda^{\abs{I \cap \c E}} \leq \( (1 + \mu)^d + (1 + \lambda)^d - 1\)^{N/2d}
\end{equation}
for all $\lambda, \mu \geq 1$. This was later relaxed to $\lambda, \mu \geq 0$ by Galvin and Tetali \cite{GT}, who studied the general problem of counting weighted graph homomorphisms, again using entropy methods. We will say more about this in a moment. See \cite{MT} for more applications of entropy to counting independent sets and graph homomorphisms.

Although the entropy method is powerful, it would be nice to have a more elementary and combinatorial proof of the bipartite case. As Kahn \cite{Kahn} writes, ``one would think that this simple and natural conjecture \dots would have a simple and natural proof.''

Here we offer some ideas for a combinatorial proof of the bipartite case that is in the same spirit as the proof presented in this paper. For a nonnegative integer $t$, let $t \cdot G$ denote a disjoint union of $t$ copies of $G$. Then $P(\lambda, t \cdot G) = P(\lambda, G)^t$. So it would suffice to prove that
\begin{equation} \label{eq:expand}
	P(\lambda, d \cdot G) \leq P\( \lambda, \frac{N}{2} \cdot K_{d,d} \).
\end{equation}
We conjecture that \eqref{eq:expand} is true term-by-term, that is, for every $k$, the number of the independent sets of size $k$ in $d \cdot G$ is at most that of $\frac{N}{2} \cdot K_{d,d}$. A slightly stronger form of this conjecture has been stated by Kahn \cite{Kahn}. Asymptotic evidence for Kahn's conjecture was later provided by Carroll, Galvin, and Tetali \cite{CGT}.

We note that the term-by-term conjecture is false for the bi-weighted version in the sense of \eqref{eq:biweighted}. Indeed, for the 2-regular bipartite graph $G = K_3 \x K_2$, there exists two independent sets in $2 \cdot G$ with 3 vertices on each side of the bipartition, but there do not exist such independent sets in $3 \cdot K_{2,2}$.

%

\subsection*{Graph homomorphisms from non-bipartite graphs.}

Galvin and Tetali \cite{GT} generalized Kahn's result and showed that for any $d$-regular, $N$-vertex bipartite graph $G$, and any graph $H$ (possibly with self-loops),
\begin{equation} \label{eq:hom}
	\abs{\Hom(G, H)} \leq \abs{\Hom(K_{d,d}, H)}^{N/2d},
\end{equation}
where $\Hom(G, H)$ denotes the set of graph homomorphism from $G$ to $H$, i.e., maps $f : V(G) \to V(H)$ where $f(u)f(v)$ is an edge of $H$ whenever $uv$ is an edge of $G$. Graph homomorphisms generalize the notion of independent sets, since we can take $H$ to be the graph $K_2$ with a self-loop adjoined to one vertex. 
In fact \cite{GT} gives a more general weighted version of \eqref{eq:hom}. Here each $i \in V(H)$ is assigned some ``activity'' $\lambda_i$. Write $\Lambda$ for the vector of activities. For each $f \in \Hom(G, H)$, let
\[
	w^\Lambda(f) = \prod_{v \in V(G)} \lambda_{f(v)},
\]
the \emph{weight} of $f$, and set
\[
	Z^\Lambda(G, H) = \sum_{f \in \Hom(G, H)} w^\Lambda(f).
\]
(The motivation for this setup comes from statistical mechanics.) It was shown in \cite{GT} that for any $d$-regular, $N$-vertex bipartite $G$, any graph $H$, and any system $\Lambda$ of positive activities on $V(H)$, we have, generalizing \eqref{eq:hom}
\begin{equation} \label{eq:hom-weighted}
	Z^\Lambda(G, H) \leq \(Z^\Lambda(K_{d,d}, H)\)^{N/2d}.
\end{equation}
A biweighted version generalizing \eqref{eq:biweighted} was also given in \cite{GT}. It was conjectured that \eqref{eq:hom-weighted} holds also when $G$ is non-bipartite. 
It is natural to try to generalize our approach in Section \ref{sec:proof} to prove this, but one runs into the problem that, unlike in the case of independent sets, there is no weight-preserving injection from $\Hom(2 \cdot G, H)$ to $\Hom(G \x K_2, H)$. For example, take $G = H = K_3$. Then there are $36$ elements in $\Hom(2 \cdot K_3, K_3)$ with weight $\lambda_1^2\lambda_2^2\lambda_3^2$ whereas there are only $24$ elements in $\Hom(K_3 \x K_2, H)$ with that weight. However, the silver lining is that there does exist a weight-preserving injection from $\Hom(2\ell \cdot K_3, K_3)$ to $\Hom(\ell \cdot K_3 \x K_2, K_3)$ for all $\ell \geq 5$. So there is hope that our method could still work.

\section*{Acknowledgements}
This research was carried out at the University of Minnesota Duluth under the supervision of Joseph Gallian with the financial support
of the National Science Foundation (grant number DMS 0754106), the National Security Agency (grant number H98230-06-1-0013), and the MIT Department of Mathematics. The author would like to thank Joseph Gallian for his encouragement and support. The author would also like to thank David Galvin for helpful email discussions, and Reid Barton, Jeff Kahn, and the anonymous referee for valuable suggestions on the exposition.

\bibliographystyle{amsplain}
\bibliography{../references}

\providecommand{\bysame}{\leavevmode\hbox to3em{\hrulefill}\thinspace}
\providecommand{\MR}{\relax\ifhmode\unskip\space\fi MR }
\providecommand{\MRhref}[2]{%
  \href{http://www.ams.org/mathscinet-getitem?mr=#1}{#2}
}
\providecommand{\href}[2]{#2}
\begin{thebibliography}{10}

\bibitem{Alon91}
N.~Alon, \emph{Independent sets in regular graphs and sum-free subsets of
  finite groups}, Israel J. Math. \textbf{73} (1991), no.~2, 247--256.

\bibitem{BergSteif}
J.~van~den Berg and J.~E. Steif, \emph{Percolation and the hardcore lattice gas
  model}, Stoch. Proc. Appl. \textbf{59} (1994), 179--197.

\bibitem{CE}
P.~J. Cameron and P.~Erd{\H{o}}s, \emph{On the number of sets of integers with
  various properties}, Number theory ({B}anff, {AB}, 1988), de Gruyter, Berlin,
  1990, pp.~61--79.

\bibitem{CGT}
T.~Carroll, D.~Galvin, and P.~Tetali, \emph{Matchings and independent sets of a
  fixed size in regular graphs}, J. Combin. Theory Ser. A \textbf{116} (2009),
  1219--1227.

\bibitem{GalvinPC}
D.~Galvin, \emph{personal communication}.

\bibitem{Galvin}
D.~Galvin, \emph{An upper bound for the number of independent sets in regular
  graphs}, Discrete Math. (to appear).

\bibitem{GT}
D.~Galvin and P.~Tetali, \emph{On weighted graph homomorphisms}, Graphs,
  morphisms and statistical physics, DIMACS Ser. Discrete Math. Theoret.
  Comput. Sci., vol.~63, Amer. Math. Soc., Providence, RI, 2004, pp.~97--104.

\bibitem{GZ}
D.~Galvin and Y.~Zhao, \emph{The number of independent sets in a graph with
  small maximum degree}, submitted.

\bibitem{CE-Green}
B.~Green, \emph{The {C}ameron-{E}rd{\H o}s conjecture}, Bull. London Math. Soc.
  \textbf{36} (2004), no.~6, 769--778.

\bibitem{GutHar}
I.~Gutman and F.~Harary, \emph{Generalizations of the matching polynomial},
  Utilitas Math. \textbf{24} (1983), 97--106.

\bibitem{Haggstrom}
O.~H\"{a}ggstr\"{o}m, \emph{Ergodicity of the hard-core model on $\mathbb{Z}^2$
  with parity-dependent activities}, Ark. Mat. \textbf{35} (1997), 171--184.

\bibitem{Kahn}
J.~Kahn, \emph{An entropy approach to the hard-core model on bipartite graphs},
  Combin. Probab. Comput. \textbf{10} (2001), no.~3, 219--237.

\bibitem{Kahn2}
\bysame, \emph{Entropy, independent sets and antichains: a new approach to
  {D}edekind's problem}, Proc. Amer. Math. Soc. \textbf{130} (2002), no.~2,
  371--378 (electronic).

\bibitem{MT}
M.~Madiman and P.~Tetali, \emph{Information inequalities for joint
  distributions, with interpretations and applications}, IEEE Trans. on
  Information Theory (to appear).

\bibitem{Sapo}
A.~A. Sapozhenko, \emph{On the number of independent sets in extenders},
  Diskret. Mat. \textbf{13} (2001), no.~1, 56--62.

\bibitem{CE-Sap}
\bysame, \emph{The {C}ameron-{E}rd{\H o}s conjecture}, Dokl. Akad. Nauk
  \textbf{393} (2003), no.~6, 749--752.

\end{thebibliography}

\end{document}